\documentclass{article} 

\usepackage{graphicx}
\usepackage{amsmath}
\usepackage{amsfonts}
\usepackage{amssymb}

\newcommand{\be}[1]{\begin{equation}\label{#1}}
\newcommand{\ee}{\end{equation}}
\newtheorem{theorem}{Theorem}[section]
\newtheorem{corollary}{Corollary}
\newtheorem{lemma}[theorem]{Lemma}

\begin{document}

\title{Dynamical system modeling fermionic limit}

\author{Dorota Bors and Robert Sta\'nczy}



\maketitle

\centerline{\scshape Dorota Bors}
\medskip
{\footnotesize
 \centerline{Faculty of Mathematics and Computer Science}
   \centerline{University of Lodz}
   \centerline{Banacha 22}
   \centerline{90-238 \L\'od\'z, Poland}
} 

\medskip

\centerline{\scshape Robert Sta\'nczy}
\medskip
{\footnotesize
 \centerline{ Instytut Matematyczny}
 \centerline{Uniwersytet Wroc{\l}awski}
 \centerline{pl. Grunwaldzki 2/4}
 \centerline{50--384 Wroc{\l}aw, Poland}
}


\thanks{}

\bigskip

\begin{abstract}

The existence of multiple radial solutions to the elliptic equation modeling
fermionic cloud of interacting particles is proved for the limiting Planck constant and intermediate values of mass parameters. It is achieved by considering the related nonautonomous dynamical system for which the passage to the limit can be established due to the continuity of the solutions with respect to the parameter going to zero.
\end{abstract}

\section{Introduction and motivation}

Consider the following elliptic boundary value problem
\begin{equation}\label{ell}
\Delta \phi (u) = H_\eta^{-1}(c-\phi (u))
\end{equation}
where $\phi$ plays the role of the the gravitational potential generated by the cloud of diffusive particles with the self--agreed density $H_\eta^{-1}(\phi (u)+c)$ distributed over $u\in B(0,1)\subset {\mathbb R}^d$ and  the constant $c$  satisfying, for given $M>0$, the mass constraint 
$$
\int_{B(0,1)} H_\eta^{-1}(c-\phi(u))\,du \,=\,M\,.
$$
The origins of the function $H_\eta$ stems from the statistical mechanics approach. The function $H_\eta$ is given and depends on the parameter $\eta \ge 0$. The form of $H_\eta$ encompasses the models arising from the Maxwell--Boltzmann and the Fermi--Dirac statistics.

We shall prove the multiplicity results for the above nonlocal BVP for the intermediate values of the mass parameter $M>0$ while the parameter $\eta>0$ is taken sufficiently close to zero. Thus it can be seen as a singular perturbation of the Maxwell--Boltzmann statistics with $\eta=0$.

The problem can be reduced, by appropriate substitution, to some dynamical system stated in (\ref{Syst:Nonau}) for given $H_\eta$ by defining the new nonlinearity $R_\eta$ as
$$
H'_\eta(z)R_\eta(z)=1\,.
$$

We consider the following functions originating from the statistical mechanics:
\begin{itemize}
\item $R_0(z)=z$ in the Maxwell--Boltzmann model with $H_0(z)=\log (z)$,
\item $R_\eta(z)=(1/z+\eta/z^{1/d})^{-1}$ in the simplified Fermi--Dirac model with $$H_\eta(z)=\log (z)+\eta z^{1-1/d}$$,
\item $R_\eta(z)=\frac{\mu(d-2)}{4}f_{d/2-2}(f_{d/2-1}^{-1}(2z/\mu))$ in the Fermi--Dirac model with 
$$\eta \mu^{2/d}=2d^{2/d-1}$$ 
and the Fermi functions $f_\alpha$ defined as
$$f_\alpha(z)=\int_0^{\infty}\frac{x^\alpha}{1+\exp (x-z)} dx\,. $$
\end{itemize}

The Fermi--Dirac model was introduced to describe in a better way the existence of the galaxies or the gaseous stars than the Maxwell--Boltzmann model. In the Maxwell--Boltzmann model the existence of blowing--up solutions for (\ref{evo} -- the so called the gravo-thermal castastrophe was proved. It was accompanied by the lack of steady states for massive clouds but was not supported by observations of evolving galaxies or stars towards stable steady states, cf. \cite{R}. The motivation for considering such form of equations comes from the models of self--gravitating diffusive particles introduced by Chavanis et al. in \cite{CSR} and developed further in \cite{C, MR2092680}.

Relating the potential $\phi$ to the new variables $x$ reduced mass and  $y$ the energy leads to the possibly nonautonomous system
\be{Syst:Nonau}\left\{\begin{array}{l}
x'(s)=(2-d)\,x(s)+\,y(s)\,,\\[6pt]
y'(s)=2\,y(s)-\,x(s)\,e^{2s}\,R_\eta(\,e^{-2s}\,y(s))\,,
\end{array}\right.\ee
with parameters $d\in {\mathbb N}\cap [3,9], \eta\ge 0$ that reduces for $R_0={\rm I}$ to the autonomous one
\be{Syst:Auton}\left\{\begin{array}{l}
x'(s)=(2-d)\,x(s)+\,y(s)\,,\\[6pt]
y'(s)=(2-\,x(s))\,y(s)\,.
\end{array}\right.\ee

Indeed the system (\ref{Syst:Nonau}) can be derived from the elliptic equation, up to constant studied in \cite{SRS, SAA}, by considering
\be{ODE}
-\,Q''+(d-1)\,r^{-1}\,\,Q'=Q\,R_\eta(r^{1-d}\,Q')
\ee
with $Q(0)=0, Q(1)=\sigma_d^{-1}M$ using the substitution relating $s, x, y$ to $r, Q, Q'$ given by
$$
Q(e^s)=x(s)\,e^{(d-2)s}, Q'(e^s)=y(s)\,e^{(d-3)s}\,.
$$
The latter equation (\ref{ODE}) describes 
$$Q(r)=\sigma_d^{-1}\int_{B(0,r)}\rho (u) du$$
the averaged (differing thus by a constant $\sigma_d$ measure of the unit sphere from notation adopted in \cite{SRS} and \cite{SAA}), i.e. integrated over the ball $B(0,r)$, the radial density $$\rho (u)=H_\eta^{-1}(c-\phi (u))$$ of the particles preserving mass $\int_{B(0,1)}\rho (u) du=M$. 

The $x$ variable is related to the rescaled mass parameter, while $y$ can be vaguely referred to the energy of the system. The precise reference is stated in the sequel.  One should note that while the system (\ref{Syst:Auton}) referred to as Maxwell--Boltzmann case is well understood as been thoroughly examined in many papers, cf. \cite{BHN, KN} and references therein, the so called Fermi--Dirac like system (\ref{Syst:Nonau}) is less studied and not many results are available, cf. \cite{SRS, SAA}. This difficulty is generated by the nonlinear nature of the $R_\eta$ function causing some additional problems and posing some extra difficulties. The problem can be also studied in slightly more general framework allowing $R_\eta$ satisfying some condition cf. Theorem 3.1 encompassing also the Fermi--Dirac case. It should be noted that the results obtained for both models differ significantly for $d=3$ and large values of mass parameter, while for small and intermediate values of mass parameter they share the common features provided the parameter $\eta$ related to the Planck constant is small enough. The main result of this paper is the convergence of properly chosen solutions of the system (\ref{Syst:Nonau}) towards the solutions to (\ref{Syst:Auton}) as $\eta \rightarrow 0$ and the mass parameter $M$ attains some intermediate values. This results in the existence of multiple solutions for the Fermi--Dirac model for $\eta$ small enough and properly chosen mass parameter $M>0$ with intermediate values as in the Maxwell--Boltzmann case. This can be depicted in the phase diagram on Figure 1 illustrating the main Theorem 3.2 and Corollary 1 of the manuscript. The results presented in this paper can be seen as continuous dependence of the solutions to the dynamical system on the parameter $\eta \ge 0$ but only for sufficiently small values of the parameter. It should be underlined that solutions for the dynamical systems are defined on the non-compact interval. Moreover, we choose some special family of the solutions characterized by the limit at minus infinity, not the whole set of possible solutions. The continuous dependence on parameters of the whole set of solutions for elliptic equations was established among others in \cite{B, BWN, BWC, BW}. One should point out that the passage to the singular limit was rigorously verified both for the related Navier--Stokes--Fourier--Poisson system by Lauren\c{c}ot and Feireisl in \cite{FL} while Golse and Saint--Raymond in \cite{GSR} dealt with celebrated Navier--Stokes and Boltzmann equations.

The solutions of the BVP with elliptic equation considered above \ref{ell} can be seen as steady states for the evolutions of the potential of particles with the density $\rho$ and  with no flux boundary condition evolving by
\begin{equation}\label{evo}
\rho_t=\nabla \cdot N\left( \theta P_\eta' \nabla \rho + \rho\nabla \Delta^{-1} \rho \right)\,,  
\end{equation}
with some positive coefficient $N$ possibly depending on other variables, where
$$
P_\eta'(z)=H'_\eta(z)z\,.
$$

\section{Derivation of the dynamical systems}

The results are the extension of the results obtained for the case $d=3$ in \cite{DSD} to higher dimension $3\le d\le 9$ and more general pressure formulae $P_\eta$ generating via $$P_\eta'(z)=zH_\eta'(z),\, H_\eta'(z)R_\eta(z)=1$$ with the function $R_\eta$ appearing in the system (\ref{Syst:Nonau}) while $z=\rho\theta^{-d/2}$ where $\theta$ is the temperature of the system and $\eta$ is the parameter related to the Planck constant.

Let us analyze the limit system (\ref{Syst:Auton}) for which the point $(0,0)$ is a saddle, while the other stationary point $(2,2(d-2))$ can change character if any $d$ is considered but if $3\le d \le 9$ then it is a sink and a Lyapunov function
$$
L(x,y)=\frac{1}{2}(x-2)^2+y-2(d-2)-2(d-2)\log(y/(2d-4))\,
$$
governs convergence towards this point as was established in \cite{BHN} and started in \cite{KN}. Indeed multiplying 
the equations (\ref{Syst:Auton}) for $x'$ by $x-2$ and $y'/y$ by $2(2-d)$ and summing them with added $y'$ one obtains
$$
\frac{d}{dt}L(x(t),y(t))=x'(t)(x(t)-2)+y'(t)-2(d-2)y'(t)/y(t)=-(x(t)-2)^2\le 0\,.
$$
Moreover, using Taylor expansion in the neighborhood of $(x,y)\sim (2,2(d-2))$ we can see that
$$L(x,y)\sim \frac12 (x-2)^2+\frac1{4(d-2)}(y-2(d-2))^2\,.$$

Furthermore, note that the condition $Q(0)=0$ can be translated to
$$
\lim_{s\rightarrow -\infty} x(s)e^{(d-2)s} = 0\,,
$$
while assuming $\rho\in L^{\infty}$ guarantees $Q(r)$ be of order $r^{d}$ at zero thus assuring $x(s)e^{-2s}$ to be bounded. Moreover, if $\rho$ is continuous then the following limit exists and is finite
$$
\lim_{s\rightarrow -\infty} x(s)e^{-2s}<\infty\,.
$$
Additionally,
$$
\rho(0)=|\rho|_\infty=\lim_{s\rightarrow -\infty} y(s)e^{-2s}<\infty\,.
$$
One assumes $R_\eta$ to be continuous on $[0,\infty)$ to claim the following lemma in the first, positive quadrant.
\begin{lemma}\label{lem}
For any solution $(x,y)$ to (\ref{Syst:Nonau}), finite $\rho_0=\lim_{s\rightarrow -\infty} y(s)e^{-2s}$ implies $$\lim_{s\rightarrow -\infty}\frac{x(s)}{y(s)}=\frac1d\,.$$
\end{lemma} 
{\bf Proof.} Using de l'Hospital rule together with the system (\ref{Syst:Nonau}) one gets the claim by
$$M=\lim_{s\rightarrow -\infty}\frac{x(s)}{y(s)}=\lim_{s\rightarrow -\infty}\frac{x'(s)}{y'(s)}=\lim_{s\rightarrow -\infty}\frac{1+(2-d)\frac{x(s)}{y(s)}}{2-e^{2s}R_\eta(e^{-2s}y(s))\frac{x(s)}{y(s)}}=\frac{1+(2-d)M}{2}$$

\section{Convergence and multiplicity results}

Consider the system describing the evolution of the difference 
\be{Syst:Deffi}\left\{\begin{array}{l}
w_\eta=x_\eta-x_0\\[6pt]
v_\eta=y_\eta-y_0
\end{array}\right.\ee
of solutions $(x_\eta,y_\eta)$ to (\ref{Syst:Nonau}) and  $(x_0,y_0)$ to (\ref{Syst:Auton}) of the form
\be{Syst:Differ}\left\{\begin{array}{l}
w_\eta'=(2-d)\,w_\eta+\,v_\eta\\[6pt]
v_\eta'=(2-x_0)\,v_\eta-y_\eta w_\eta-\,x_\eta\,e^{2s}\,S_\eta(\,e^{-2s}\,y_\eta)
\end{array}\right.\ee
where
\begin{equation}
S_\eta(z)=z-R_\eta(z)\,.
\end{equation}

Now we shall prove crucial a priori bound for the term $x_\eta\,e^{2s}\,S_\eta(\,e^{-2s}\,y_\eta)$ appearing in (\ref{Syst:Differ}). Set $\rho_0>0$ and take for any $\rho\le \rho_0$ the solution $y$ such that 
$$
\rho=\lim_{s\rightarrow -\infty} y(s)e^{-2s}\,.
$$
Then by Lemma \ref{lem} we have that $y_\eta e^{-2s}\nearrow \rho$ and $x_\eta e^{-2s}\nearrow \frac1d\rho$ as $s\rightarrow -\infty$ hence
$y_\eta \le \rho_0e^{2s}$ and $x_\eta \le \frac1d\rho_0 e^{2s}$ whence
\begin{equation}
dx_\eta\,e^{2s}\,S_\eta(\,e^{-2s}\,y_\eta) \le \rho_0 e^{4s} \max_{[0,\rho_0]} S_\eta =\rho_0 e^{4s} \overline{S_{\eta,\rho_0}}= \rho_0 \overline{S_{\eta,\rho_0}}\,,
\end{equation}
where $\overline{S_{\eta,\rho_0}}$ is increasing in $\rho_0$ and decreasing to zero as $\eta$ tends to $0$.

Multiplying $w'_\eta$ by $w_\eta$ and $v'_\eta$ by $v_\eta$ respectively one obtains
$$
w'_\eta w_\eta=(2-d)w_\eta^2+v_\eta w_\eta
$$
and
$$
v'_\eta v_\eta=(2-x_0)v_\eta^2-y_\eta w_\eta v_\eta - x_\eta e^{2s}S_\eta(e^{-2s}y_\eta) v_\eta\,.
$$
Next setting $\chi = w_\eta^2 + v_\eta^2$ one obtains
$$ \chi' \le \alpha \chi + \beta$$
where $\alpha$ and $\beta$ are the terms bounded with respect to $\eta$. Indeed one can estimate
$$
|(1-y_\eta)v_\eta w_\eta| \le \frac12 (w_\eta^2+v_\eta^2) \max \{1,\max y_\eta -1\}
$$
and
$$
|2-x_0|\le \max \{2,\max x_0-2 \}
$$
while
$$
y_\eta \le \rho_0, x_\eta\le \frac1d\rho_0, |v_\eta| \le \frac12 (1+v_\eta^2)\,.
$$
Thus coming back to the estimate on $\chi$ and by the Gronwall lemma
$$\chi \le \beta e^\alpha\,.$$
But it should be noted that 
$$\beta \le \frac12 x_\eta e^{2s}S_\eta\le \frac{1}{2d} \rho_0 S_\eta\le \frac{1}{2d} \rho_0 \overline{S_{\eta,\rho_0}}$$
while
$S_\eta(z) \le C(\eta)D(z)$
where $C(\eta)\rightarrow 0$ as $\eta\rightarrow 0$.

One can alternatively proceed with the same conclusion as in \cite{DSD} defining $A_\eta=\sup_{s\le t}e^{-2s}|x_\eta(s)-x_0(s)|\,, B_\eta=\sup_{s\le t}e^{-2s}|y_\eta(s)-y_0(s)|$ getting almost everywhere
$$ \frac{d}{dt}(e^{dt}A_\eta) \le e^{dt} B_\eta$$
and with some constant $\kappa$ one gets
$$ \frac{d}{dt} B_\eta \le (\rho_0A_\eta+\frac{1}{d}\rho_0B_\eta+\eta \kappa)\,.$$
Then integration over $(-\infty,t)$ yields
$$dA_\eta \le B_\eta \le \frac12 (\rho_0A_\eta+\frac{1}{d}\rho_0 B_\eta + \eta \kappa)e^{2t}.$$
Finally, as before, using a Gronwall estimate, one gets as required
$$0 \le d A_\eta \le B_\eta\le \frac12 \eta \kappa e^{\rho_0/d}\,.$$

Thus we have proved the following convergence theorem
\begin{theorem}\label{conv}
Fix any natural number $3\le d\le 9$ and $\rho_0>0$ and take any $\rho \in (0,\rho_0]$ such that $\lim_{s\rightarrow\-\infty} y(s)e^{-2s}=\rho$ for some solution to the system. Assume that the continuous function $R_\eta$ satisfies
$$
0\le z-R_\eta(z)\le C(\eta)D(z)
$$
where $C(\eta)\rightarrow 0$ as $\eta\rightarrow 0$ and $D:[0,\infty)\rightarrow [0,\infty)$ is a continuous function such that $D(0)=0$. Then the solution $(x_\eta, y_\eta)$ converges uniformly to $(x_0,y_0)$ on $(-\infty,0]$ and in particular $x_\eta(0)$ converges to $x_0(0)$.
\end{theorem}

Recall from \cite{DSD} in the case $d=3$ the generic for $3\le d \le 9$ phase portrait for the Maxwell--Boltzmann case with $R_0= I$ identity function and $(2,2(d-2))=(2,2)$.

\begin{figure}[ht]
\label{Fig:Max}
\begin{center}
\includegraphics[height=3.5cm]{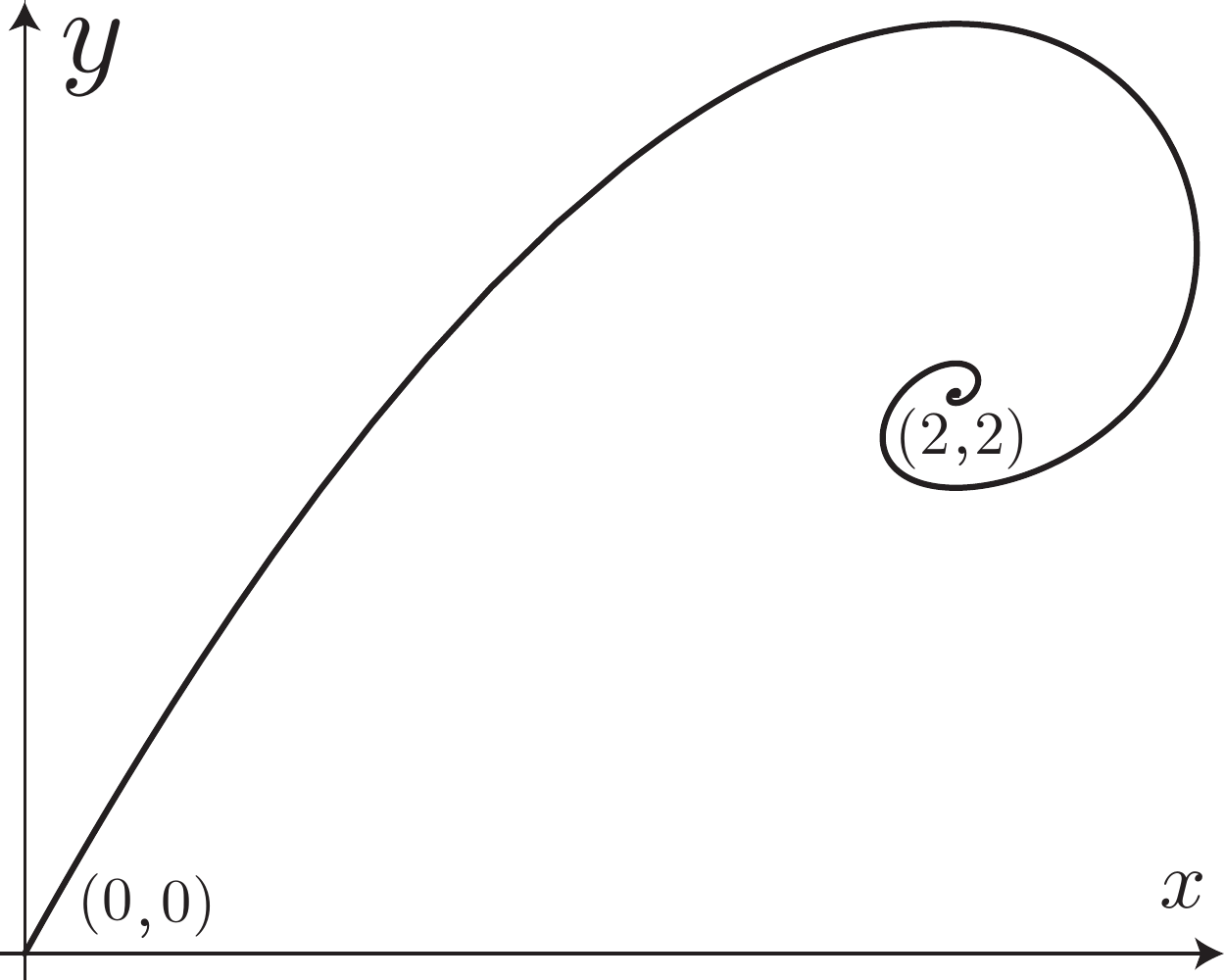}\hspace*{30pt}\includegraphics[height=3.5cm]{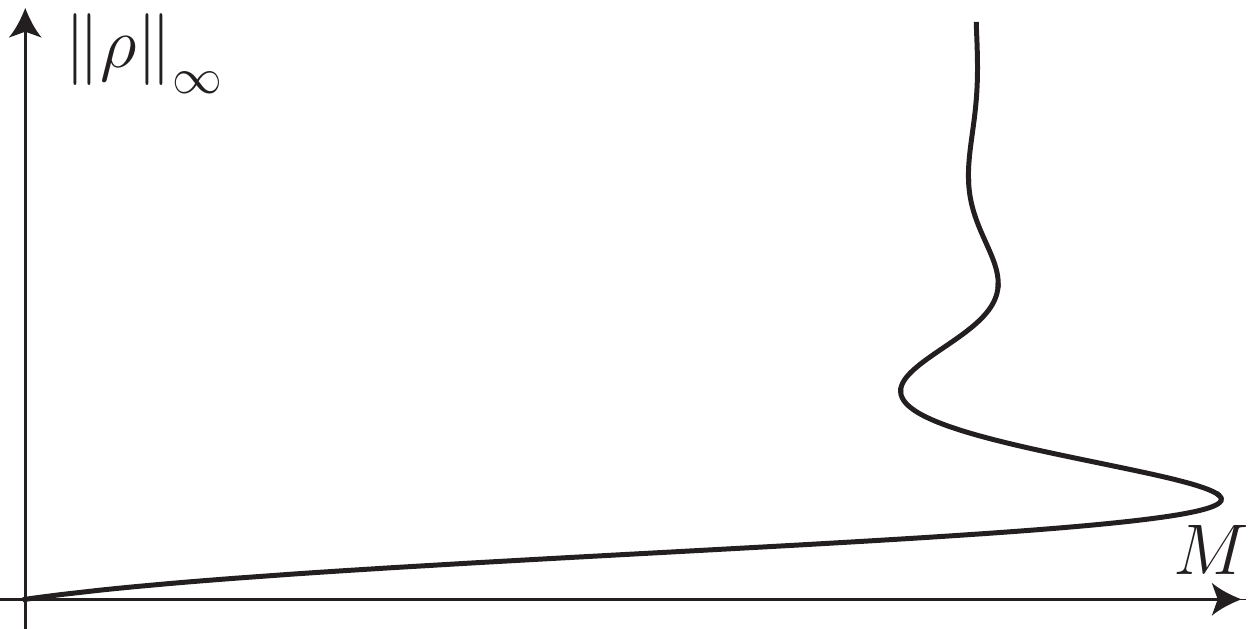}
\end{center}
\caption{\it\small Left: the heteroclinic orbit joining the points $(0,0)$ and $(2,2)$ in the Maxwell--Boltzmann case. Right: the mass--density diagram.}
\end{figure}

An easy corollary, due to the fact that the mass of the system is related to $x_\eta(0)$ of the Theorem \ref{conv}, can be formulated as follows.
 
\begin{theorem}
For any mass parameter $M$ in the corresponding, intermediate range for Fermi--Dirac like models modeled by $R_\eta$, with $\eta>0$ small enough, satisfying the condition from the Theorem \ref{conv} there exists as many solutions as for the Maxwell--Boltzmann case with $R_0=I$ depicted in the Figure 1 and depending on the intersection of the vertical line (setting thus the mass $M>0$) with the bifurcation curves.
\end{theorem}

The details of the proof are the same as in \cite{DSD} and are omitted herein but they focus on the continuous dependence of the mass $M$ on the density $\rho$ or in other words $x_\eta(0)$ on $\lim_{s\rightarrow -\infty}e^{-2s}y(s)$ expressed in the language of the dynamical system variables.

\begin{corollary}
For the intermediate values of the mass parameter $M$ there exists multiple solutions to the Fermi--Dirac $R_\eta$ and generalizations obeying the condition from Theorem \ref{conv} provided the $\eta$ parameter is sufficiently small.
\end{corollary}

We show that the phenomena appearing in Fermi--Dirac model for large values of the mass (cf. \cite{SDI}) parameter differentiating between dimensions $d=3$ (solution for any mass), $d\ge 5$ (solution only up to some mass parameter) are not present for intermediate value of mass parameter, where the behavior is generic for any dimension $3\le d\le 9$. This is accompanied by existence of multiple solutions for any dimension if we are close enough to the Maxwell--Boltzmann case with $\eta=0$. One should note also that for small values of mass parameter the uniqueness holds as was noted in \cite{DOLBEAULT:2009:HAL-00349574:2}.

\section{Appendix}

\begin{lemma}
For simplified Fermi--Dirac model we have straightforward estimate
$$
z^{-1-2/d}(z-R_\eta(z))\le 1\,.
$$
\end{lemma}

\begin{lemma}
For the Fermi--Dirac model we have
$$
z^{-1-2/d}(z-R_\eta(z))\le \left(\frac{2}{\mu}\right)^{2/d}C\,,
$$
where
$$
C=\max_{w\in [0,\infty)} w^{-1-2/d}(w-\frac{d-2}{2}\zeta (w))
$$
and
$$
\zeta(w) =f_{d/2-2}(f_{d/2-1}^{-1}(w))\,.
$$
\end{lemma}
{\bf Proof.} Notice that due to the asymptotics of the Fermi functions \cite{BNS}
$$\frac{z-R_\eta(z)}{z^{1+2/d}}\cdot\left(\frac{\mu}{2}\right)^{1+2/d}=\frac{\frac{\mu}{2}\cdot\frac{2z}{\mu}-\frac{\mu(d-2)}{4}\zeta(2z/\mu)}{\left(\frac{2z}{\mu}\right)^{1+2/d}}=\frac{\mu}{2}\frac{w-\frac{d-2}{2}\zeta(w)}{w^{1+2/d}}\le \frac{\mu}{2}C\,.$$
Recall the relation between constants that appear above to agree behavior of the functions at $\infty$
$$
\eta \mu^{2/d}=2d^{2/d-1}\,.
$$
Hence $\mu\rightarrow \infty$ when $\eta\rightarrow 0^{+}\,.$
Indeed using the estimates from \cite{BNS} or \cite{SM} we have that
$$
f_\alpha(w) \sim \frac1{\alpha+1} w^{\alpha+1}, w\sim \infty
$$
while
$$
f_\alpha(w) \sim  \Gamma (\alpha+1) \exp(w) , w\sim 0^{+}\,.
$$

The motivation for considering the pressure $p$ in the model equation
\begin{equation}
\rho_t=\nabla \cdot N\left( \nabla p + \rho\nabla \Delta^{-1} \rho \right)\,,  
\end{equation}
in the form (\ref{evo}) with the specific dependence on the temperature $\theta$, the density $\rho$ and the dimension of the ambient space $d$ reading
$$p(\theta,\rho)=\theta^{d/2+1}P(\rho\theta^{-d/2})$$
with some given $P$ function (we drop dependence on $\eta$) is threefold. First of all one can for $zH'(z)=P'(z)$ with $z=\rho\theta^{-d/2}$ establish the entropy formula 
$$
{\mathcal W}=\int_{B(0,1)} \left(\rho H(\rho\theta^{-d/2}) - \left(\frac{d}{2}+1\right)\theta^{d/2}P(\rho\theta^{-d/2})\right)
$$
due to this assumption on the pressure form, cf. \cite{BS}. Then the number of astrophysically motivated examples can be found as: Maxwell--Boltzmann, Fermi--Dirac, Bose--Einstein or polytropic statistics modeling clouds of particles, galaxies or stars. Finally, some monoatomic gases require this assumption which can be found in \cite{EGH, F, FN, MR}. To this end we recall for $d=3$ Maxwell's equation with kinetic internal energy per molecule $e$
$$\rho^2e_\rho=p-\theta p_\theta\,.$$
While for monoatomic gas the relation holds
$$3p=2\rho e\,.$$
Hence
$$3p_\rho=2e+2\rho e_\rho$$ plugged into Maxwell's equation derived from the Gibb's relation (cf. \cite{F}) yields 
$$2p-2\theta p_\theta=2\rho^2e_\rho=3\rho p_\rho-2e\rho=3\rho p_\rho-3p\,.$$
This gives the linear first order partial differential equation 
$$5p=3\rho p_\rho+2\theta p_\theta$$
that can be solved with characteristics i.e. the system of equations
$$\rho'=3\rho, \theta'=2\theta, p'=5p$$
with two first integrals of the form
$$p\theta^{-5/2}, \rho\theta^{-3/2}\,.$$
This yields the solution in the implicit form
$$\Phi(p\theta^{-5/2},\rho\theta^{-3/2})=0$$
or explicit form
$$p\theta^{-5/2}=P(\rho\theta^{-3/2})$$
giving
$$p=\theta^{5/2}P(\rho\theta^{-3/2}).$$
In higher dimension replacing $3$ with $d$ would yield the corresponding formula
$$p=\theta^{d/2+1}P(\rho\theta^{-d/2}).$$

\section{Open problems and possible extensions}
One can consider the Dirichlet boundary value problem with elliptic equation
$$\Delta \phi (u)= \rho (u) = H_\eta^{-1}(c-\phi (u))$$
where the constant $c$ is chosen so that the mass constraint holds
$$\int_{B(0,1)} \rho (u) du = \int_{B(0,1)} H_\eta^{-1}(c-\phi(u))du = M\,.$$
The entropy can be used in dimensions $d=3$ with any mass $M>0$ or $d=4$ and the mass $M$ sufficiently small to obtain the minimizer solving the related Euler-Lagrange equation.
To be more specific the dual approach, cf. \cite{SDI}, uses the neg-entropy functional
$$
{\mathcal V}=\int_{B(0,1)} \left(\rho H_\eta(\rho \theta^{-d/2}) - \theta^{d/2}P_\eta(\rho\theta^{-d/2})+\frac{1}{2\theta}\rho \Delta^{-1}\rho \right)
$$
over the space of integrable functions $\rho \in L^{1+2/d}$. The functional is coercive and can be decomposed into compact and continuous part and lower--semicontinuous and convex part thus making the direct approach feasible to yield the existence of minimizer. It seems that the results of \cite{B, BWN, BWC, BW} can be used to get the continuity of the set of minimizers at least for sufficiently small mass $M>0$. The only obstacle is that the limiting functional is defined over the space of $\rho \log \rho$ integrable functions as $\eta\rightarrow 0^{+}$.

Moreover, one can consider with necessary modifications the following nonlinearities
\begin{itemize}
\item $R_\eta(z)=\frac{\mu(d-2)}{4}g_{d/2-2}(g_{d/2-1}^{-1}(2z/\mu))$ in the Bose--Einstein model with Bose functions $g_\alpha$ defined by
$$g_\alpha(z)=\int_0^{\infty}\frac{x^\alpha}{1-\exp (x-z)} dx $$
requiring some limits for the density, or rather the ratio $\rho/\theta^{d/2}$,
\item $R$ in classical King's model, cf. \cite{CLM}, being the intermediate between Maxwell--Botlzmann and Fermi--Dirac cases\,.
\end{itemize}

\section{Acknowledgements.}
The remark on monoatomic gases that can be found in the Appendix was pointed out by Eduard Feireisl and we would like to thank him for this comment. The impact of Jean Dolbeault who contributed to the three--dimensional case \cite{DSD} could not be overestimated so the due credit is paid to him herein.



\begin{thebibliography}{99}


\bibitem{BHN}
{\sc P. Biler, D. Hilhorst, and T. Nadzieja},
Existence and nonexistence of solutions for a model of gravitational interaction of particles, II,
{\it Coll. Math.} {\bf 67} (1994), 297--308.

\bibitem{BS}
{\sc  P. Biler and R. Sta\'nczy}, 
Parabolic--elliptic systems with general density-pressure relations,
{\it Surikaisekikenkyusho Kokyuroku}, {\bf 1405} (2004), 31--53.

\bibitem{BNS}
{\sc P.~Biler, T.~Nadzieja, and R.~Sta{\'n}czy},
Nonisothermal systems of self-attracting {F}ermi-{D}irac particles,
{\it Banach Center Publ.}, {\bf 66} (2004), 61--78.

\bibitem{B}
{\sc D. Bors}, Superlinear elliptic systems with distributed and boundary controls.
{\it Control and Cybernetics}, {\bf 34} (2005), 987--1004.

\bibitem{BWN}
{\sc D. Bors and S. Walczak}, Nonlinear elliptic systems with variable boundary data.
{\it Nonlinear Analysis: Theory, Methods and Applications}, {\bf 52} (2003), 1347--1364.

\bibitem{BWC}
{\sc D. Bors and S. Walczak}, 
Stability of nonlinear elliptic systems with distributed parameters and variable boundary data. 
{\it Journal of Computational and Applied Mathematics}, {\bf 164-165} (2004), 117--130.
  
 \bibitem{BW}       
{\sc D. Bors and S. Walczak}, 
 Optimal control of elliptic systems with distributed and boundary controls. 
 {\it Nonlinear Analysis: Theory, Methods and Applications}, {\bf 63} (2005), 1367--1376.

  \bibitem{C}
{\sc P.-H. Chavanis}, 
Phase transitions in self-gravitating systems,
{\it International Journal of Modern Physics B}, {\bf 20} (2006), 3113--3198.

\bibitem{MR2092680}
{\sc P.-H. Chavanis, P.~Lauren{\c{c}}ot, and M.~Lemou}, 
Chapman-{E}nskog derivation of the generalized {S}moluchowski equation, 
{\it Phys. A},  {\bf 341} (2004), 145--164.

\bibitem{CLM}
{\sc P.-H. Chavanis, M.~Lemou, and F. M\'ehats}
Models of dark matter halos based on statistical mechanics: The classical King model,
{\it Phys. Rev. D}, {\bf 91} (2015), 063531.

\bibitem{CSR}
{\sc P.-H. Chavanis, J.~Sommeria, and R.~Robert}, {\em Statistical mechanics of
  two-dimensional vortices and collisionless stellar systems}, Astrophys. J., {\bf 471} (1996), 385.


\bibitem{DSD} 
{\sc J. Dolbeault and R. Sta\'nczy}, 
Bifurcation diagram and multiplicity for nonlocal elliptic equations modeling gravitating systems based on Fermi--Dirac statistics, {\it Disc. Cont. Dyn. Sys. A}, {\bf 35} (2015), 139--154.

\bibitem{DOLBEAULT:2009:HAL-00349574:2}
{\sc J.~{D}olbeault and R.~{S}ta\'nczy}, 
{N}on-existence and uniqueness results for supercritical semilinear elliptic equations, 
{\it {A}nnales {H}enri {P}oincar{\'e}}, {\bf 10} (2009), 1311--1333.

\bibitem{EGH}
{\sc S. Eliezer, A. Ghatak, and H. Hora}, 
{\it An Introduction to Equations of State: Theory and Applications,} 
Cambridge University Press, Cambridge, 1986.

\bibitem{FL}
{\sc E. Feireisl and P. Lauren\c ot}, 
Non-isothermal Smoluchowski–Poisson equations as a singular limit of the Navier–Stokes–Fourier–Poisson system, {\it J. Math. Pures Appl.}, {\bf 88} (2007), 325--349.

\bibitem{F}
{\sc E. Feireisl},
{\it Mathematics of complete fluid systems},
available online:\\ www.math.cas.cz/fichier/course/filepdf/course\_pdf\_20121011171111\_35.pdf

\bibitem{FN}
{\sc E. Feireisl and A. Novotn\'y},
{\it Singular limits in thermodynamics of viscous fluids},
Birkhauser-Verlag, Basel, 2009

\bibitem{GSR}
{\sc F. Golse and L. Saint--Raymond},
The Navier--Stokes limit of the Boltzmann equation for bounded collision kernels
{\it Invent. Math.}, {\bf 155} (2004), 81--161.

\bibitem{MR544879}
{\sc B.~Gidas, W.~M. Ni, and L.~Nirenberg}, 
Symmetry and related properties via the maximum principle, 
{\it Comm. Math. Phys.}, {\bf 68} (1979), 209--243.

\bibitem{GN}
{\sc M. Grendar and R.K. Niven},
Generalized Maxwell-Boltzmann, Bose-Einstein, Fermi-Dirac and Acharya-Swamy, 
{\it ArXiv} (2008).

\bibitem{MR0340701}
{\sc D.~D. Joseph and T.~S. Lundgren}, 
Quasilinear Dirichlet problems driven by positive sources, 
{\it Arch. Rational Mech. Anal.}, {\bf 49} (1972/73), 241--269.
  
\bibitem{KN} 
{\sc A. Krzywicki and T. Nadzieja}, 
Some results concerning the Poisson--Boltzmann
equation, {\it Appl. Math.}, {\bf 21} (1991), 265--272.

\bibitem{MR}
{\sc I. M\"uller and T. Ruggieri},
{\it Extended Thermodynamics}, Springer Tracts in Natural Philosophy {\bf 37}, Springer Verlag 1993.

\bibitem{R}
{\sc R. Robert}, 
On the gravitational collapse of stellar systems, 
{\it Class. Quantum Grav.}, {\bf 15} (1998) 3827--3840

\bibitem{SDI}
{\sc R. Sta\'nczy}, 
Steady states for a system describing self-gravitating
Fermi–Dirac particles, {\it Diff. Int. Equations}, {\bf 18} (2005), 567--582.

\bibitem{SRS} 
{\sc R. Sta\'nczy}, 
The existence of equlibria of many-particle systems,
{\it Proc. Roy. Soc. Edinb. A}, {\bf 139} (2009), 623--631.

\bibitem{SM}
{\sc R. Sta\'nczy}, 
On an evolution system describing self-gravitating particles in microcanonical setting,
{\it Monatshefte f\"ur Mathematik}, {\bf 162} (2011), 197--224.

\bibitem{SAA} 
{\sc R. Sta\'nczy}, 
On stationary and radially symmetric solutions to some drift-diffusion equations with nonlocal term,  
{\it Appl. Anal.}, {\bf 95} (2016), 97--104.



\end{thebibliography}
\end{document}